\documentclass[11pt]{article}

\usepackage[english]{babel}
\usepackage{graphicx}
\usepackage{amsmath, amssymb}
\usepackage{amsthm}
\usepackage{color}
\usepackage{subcaption}
\usepackage{hyperref}
\hypersetup{
    colorlinks = true,
    citecolor = {blue},
    linkcolor=blue,}
\usepackage[title]{appendix}

\newcommand{\R}{\mathbb{R}}
\newcommand{\N}{\mathbb{N}}
\newcommand{\Z}{\mathbb{Z}}

\usepackage{color, colortbl}
\definecolor{Gray}{gray}{0.9}

\newtheorem{theorem}{Theorem}[section]
\newtheorem{remark}[theorem]{Remark}

\newtheorem{proposition}[theorem]{Proposition}

\usepackage{bm}
\allowdisplaybreaks

\usepackage[top=80pt,bottom=80pt,left=60pt,right=60pt]{geometry}
\numberwithin{equation}{section}

\begin{document}



\title{Wave scattering in 1D: D'Alembert-type representations and a reconstruction method}

\author{
  Konstantinos Kalimeris\thanks{Mathematics Research Center, Academy of Athens, Greece} \,  
   and
  Leonidas Mindrinos\thanks
  {Department of Natural Resources Development and Agricultural Engineering, Agricultural University of Athens, Greece}}

\maketitle

\begin{abstract}
We derive the extension of the classical d'Alembert formula for the wave equation, which provides the analytical solution for the direct scattering problem for a medium with  constant refractive index; this is achieved by employing results obtained via the Fokas method. This methodology is further extended to a medium with piecewise constant refractive index, providing the apparatus for the solution of the associated inverse scattering problem. Hence, we provide an exact reconstruction method which is valid for both full and phaseless data.
\end{abstract}

\section{Introduction}\label{sec_intro}

We consider the scattering problem of a wave by a layered medium in the positive half line. The medium is supported in $(0,L),$ where $L>0$ describes its length. We examine both the direct and the inverse problem. The direct problem is to compute the solution of the wave equation at a detector position $x= L+D,$ with $D>0,$ given the optical properties of the medium (refractive index $\nu$ and length of every layer) and the initial wave, to be supported at $(L,\infty).$

The much more interesting, but difficult to solve, problem is the corresponding inverse problem where we aim to reconstruct the optical properties of the medium from the knowledge of the back-reflected data, meaning the wave at $x=L+D,$ for some time interval $t\in (0,T).$

The study of the inverse problem for the wave equation in layered media plays a crucial role in understanding and modelling complex wave propagation phenomena. It enables the estimation of  important material parameters in various practical applications such as seismology, geophysics, and medical imaging \cite{Bor99, Cat04, Dav19, ElbMinSch15, Som91, Sym09}.

In such cases, the media are commonly assumed to have a structure composed of multiple layers with different physical properties. For example, this assumption holds true  for the seismic layers in the Earth's subsurface but also for the human skin. The wave propagation becomes more complex due to the presence of discontinuities and variations in material properties across the different layers. 

The Fresnel equations are often relevant in the study of inverse problems for the wave equation. These equations describe the reflection and transmission of waves (acoustic and electromagnetic) at the interface between two media with different refractive indices \cite{BorWol99, hecht}.  

For instance, in optical coherence tomography, the measured intensity is the interference of a backscattered light (coming from the sample) with a back-reflected light (from a fixed mirror). For a layered medium, the backscattered signal can be described by a (infinite) sum of waves back-reflected by the discontinuities of the medium (layer interfaces) with coefficients given by the Fresnel equations. A layer stripping algorithm is usually applied for solving the inverse problem, see for example \cite{BruCha05, ElbMinVes21, Fro96, Syl96, Ves23}.

The above approach is based on an analytical representation which involves infinite series of transforms of the initial and boundary conditions. This is reflected also on the numerical simulations where one has to truncate the summation, depending on the number of multiple reflections to be considered. In the current work, we bypass this issue by deriving a d'Alembert-type representation of the solution for the wave equation on the layered medium. The d'Alembert-type representation involves only finite sums of known functions on the physical domain, which are strikingly different from the previously mentioned infinite series.

The derivation of this solution representation is based on the Fokas method for solving linear initial-boundary value problems (IBVPs). This method, which is also known as the unified transform, was introduced in 1997 by Fokas \cite{F97}, for solving IBVPs for nonlinear integrable partial differential equations (PDEs). Later, it was realised that it produced effective analytical and numerical solutions for linear PDEs \cite{F02}; for an overview of the method we refer to \cite{F08,FK22}.
The last two decades hundreds of works have been published, using and extending the Fokas method to substantially different directions, including fluid dynamics \cite{AFM06, AF11, FK17}, control theory \cite{KO20,KOD23}, regularity results \cite{Him20,OY19,BFO20}.  In general, the Fokas method yields solutions of the IBVPs as integral representations involving the Fourier (spectral) transforms of the initial and boundary conditions. Starting from these representations, one can also obtain representations in the physical space, in terms of the associated Green's functions, which involve integrals in the spectral space. For the case of the wave equation, these integrals can be easily computed explicitly, yielding, in a natural way, d'Alembert type representations. 

 The current work is based on the formulae obtained in \cite{FokKal22} where the classical d'Alembert solution of the wave equation was extended to the half-line and the finite interval problem. Employing the fact that the solutions of the IBVPs derived by the Fokas method converge uniformly to the assigned initial and boundary conditions, we extend the classical d'Alembert solution to the layered medium, by coupling the relevant conditions on the boundaries between the layers. This concept reduces the solution of the problem to a system of algebraic or functional equations, which is to be contrasted to the initial problem involving a system of partial differential equations.


The paper is organised as follows: 
In section \ref{sec_direct}, we present the characterisation of IBVP for the wave equation for the single and double-layered medium. Furthermore, we derive the extension of the d'Alembert-type solution for the single-layered medium, as well as the analogue formula for the measurement at the detector for the case that the initial condition displays a delta-like pulse distribution, both for the single and double-layered medium. Based on the above results, in section \ref{sec_inverse}, we present in detail the solution for the inverse scattering problem for the multi-layered medium,  both for full and phaseless data/measurements. We note that although the latter is an ill-posed problem, we provide evidence that one could find the unique solution if the total length of the medium is known. The numerical illustration of both the direct and inverse problem is provided in section \ref{sec_numerics}.

%

\section{The direct scattering problem}\label{sec_direct}

We consider the scattering problem in the half line $x>0.$ The support of the medium is in the interval $(0,L)$ with boundaries placed at $x=0$ and $x=L.$ We assume illumination from the right, meaning for $x>L$ with a wave having  support (short pulse) in the interval $(L,\infty).$ We define
\begin{equation}
c_j = \frac{c_0}{\nu_j}, \quad j = 1,2,\ldots,
\end{equation} 
where $c_0 \equiv 1,$ describes the speed in free space. Then, $0<c_j<1,$ is the wave speed in the $j$-th layer with refractive index $\nu_j,$ for  $j=1,2, \ldots.$

 The direct problem then reads: Given the initial wave and the properties of the medium (wave speed and length), find the back-scattered wave at some position $x>L.$ 
  
\subsection{Single-layered medium}
  
For a single-layered medium, we have only one wave speed $c_1,$ for $0<x<L,$ and the scattering problem can be decomposed into two sub-problems in the intervals $(0,\,L)$ and $(L,\,+\infty)$ as follows:

\begin{description}
\item[IBVP I (half line)]
\begin{subequations}\label{system_un}
\begin{alignat}{3}
u_{tt} - u_{xx} & = 0,  \quad &&x >L, \, t>0,  \label{eq_un1}\\
u(x,0) = U_0 (x), \, u_t (x,0) &=0,   \quad &&x >L, \label{ic_un1}\\ 
u(L,t) = G (t), \, u_x(L,t) &= F (t), \quad &&t>0, \label{bc_un1}
\end{alignat}
\end{subequations}
where $U_0 : [L,\,+\infty) \rightarrow \R.$
 
\item[IBVP II (finite interval)]
\begin{subequations}\label{system_bou}
\begin{alignat}{3}
v_{tt} - c_1^2 v_{xx} & = 0,  \quad && 0 <x <L, \, t>0,  \label{eq_b1}\\
v(x,0) = v_t (x,0) &=0,   \quad && 0 <x <L, \label{ic_b1}\\ 
v(0,t) &= 0, \quad && t>0, \label{bc_b1}\\ 
 v(L,t) = H (t), \, v_x(L,t) &= Q (t), \quad &&t>0. \label{bc_b2}
\end{alignat}
\end{subequations}
\end{description}

In addition, we impose the compatibility condition 
\begin{equation}\label{compatibility}
U_0 (L) = G(0) = H(0)=0.
\end{equation}
The problems \eqref{system_un} and  \eqref{system_bou} are coupled through the boundary conditions \eqref{bc_un1} and \eqref{bc_b2} to be seen as continuity equations in the sense:
\begin{subequations}\label{cont_eq0}
\begin{alignat}{2}
G (t) &= u(x,t) \mid_{x=L} = v(x,t) \mid_{x=L}  =  H (t) , \qquad  & t>0, \label{cont1}\\
F (t) &= u_x(x,t) \mid_{x=L} = v_x(x,t) \mid_{x=L} = Q (t) , \qquad &  t>0. \label{cont2}
\end{alignat}
\end{subequations}
 
In order to obtain the solution of \eqref{system_un} we employ \cite[Eq. (25)]{FokKal22}  using the change of variables $x= \tilde x +L:$
\begin{equation}\label{sol_un}
u(x,t) = \frac{1}{2} u_0 (x-L+t) + \begin{cases}
  \tfrac{1}{2} u_0 (x-L-t), & x>L+t, \\
  g (t-x+L) -\tfrac{1}{2} u_0 (t-x+L), & L<x<L+t.
\end{cases}
\end{equation}
This function satisfies the initial conditions
\begin{equation}
u(x,0) = u_0 (x-L), \quad u_t (x,0)=0, \quad x>L,
\end{equation}
as well as the condition $ u(L,t) =  g (t), \quad t>0$.

Hence, it suffices to set $g(t) = G(t)$ and $u_0 (x-L) = U_0 (x)$ which modifies \eqref{sol_un} into 
\begin{equation}\label{sol_un_1}
u(x,t) = \frac{1}{2} U_0 (x+t) + \begin{cases}
  \tfrac{1}{2} U_0 (x-t), & x>L+t, \\
  G (t-x+L) -\tfrac{1}{2} U_0 (t-x+2L), & L<x<L+t ,
\end{cases}
\end{equation}
which satisfies the boundary condition
\begin{equation}\label{sol_un_2}
u(L,t) =  G (t), \quad t>0.
\end{equation}
Thus, the function defined in \eqref{sol_un_1} is the solution of \eqref{system_un}.
The spatial derivative of \eqref{sol_un_1}, evaluated at $x=L,$ gives
\begin{equation}\label{sol_un_3}
u_x(L,t) =   U'_0 (L+t) -G' (t) =: F(t), \quad t>0.
\end{equation}

Similarly, from \cite[Eq. (63)]{FokKal22} setting $t = \tilde{t}/c_1$ we get
\begin{equation}\label{sol_bou}
v(x,t) = \sum_{n=0}^{\left[ \tfrac{c_1 t +x-L}{2L}\right] } h (c_1 t +x - (2n+1)L) - \sum_{n=0}^{\left[ \tfrac{c_1 t -x-L}{2L}\right] } h (c_1 t -x - (2n+1)L),
\end{equation}
where $[a]$ denotes the integer part of the real number $a$. The solution provided by the Fokas method converges uniformly to the initial and boundary conditions. Indeed, it is straightforward to observe  that the initial conditions \eqref{ic_b1} and the boundary condition \eqref{bc_b1} are satisfied. Also, for $x=L,$ we obtain 
\begin{equation}\label{sol_bou_02}
\begin{aligned}
v(L,t) &= \sum_{n=0}^{\left[ \tfrac{c_1 t}{2L}\right] } h \left(c_1 t -2nL\right) - \sum_{n=0}^{\left[ \tfrac{c_1 t}{2L}\right]-1 } h \left( c_1 t - 2 (n+1)L\right) \\
&= \sum_{n=0}^{\left[ \tfrac{c_1 t}{2L}\right] } h \left(c_1 t -2nL\right) - \sum_{n=1}^{\left[ \tfrac{c_1 t}{2L}\right] } h \left( c_1 t - 2nL\right)\\
&= h (c_1 t),
\end{aligned}
\end{equation}
where we have set $n+1 = \tilde n$ in the second sum. We set $h(c_1 t) = H(t)$ and we derive the solution of \eqref{system_bou} as
\begin{equation}\label{sol_bou_1}
v(x,t) = \sum_{n=0}^{\left[ \tfrac{c_1 t +x-L}{2L}\right] } H \left( t +\tfrac{x - (2n+1)L}{c_1}\right) - \sum_{n=0}^{\left[ \tfrac{c_1 t -x-L}{2L}\right] } H \left( t - \tfrac{x + (2n+1)L}{c_1}\right),
\end{equation}
satisfying
\begin{equation}\label{sol_bou_2}
v(L,t) = H(t), \quad t>0.
\end{equation}
 By taking the spatial derivative of \eqref{sol_bou} at $x=L,$ we obtain
\begin{equation}\label{sol_bou_3}
\begin{aligned}
v_x (L,t) &= \frac1{c_1}\sum_{n=0}^{\left[ \tfrac{c_1 t}{2L}\right] } H' \left( t -\tfrac{2nL}{c_1}\right) + \frac1{c_1}\sum_{n=0}^{\left[ \tfrac{c_1 t }{2L}\right]-1 } H' \left( t - \tfrac{2 (n+1)L}{c_1}\right) \\
&= \frac1{c_1}H' ( t) + \frac2{c_1}\sum_{n=1}^{\left[ \tfrac{c_1 t}{2L}\right] } H' \left( t  - \tfrac{2n L}{c_1} \right) =: Q(t).
\end{aligned}
\end{equation}

Then,  the continuity equations  \eqref{cont_eq0} using \eqref{sol_un_2},\eqref{sol_un_3},\eqref{sol_bou_2} and \eqref{sol_bou_3} result in the system
\begin{subequations}\label{cont_eq00}
\begin{alignat}{3}
G (t) &= H ( t), \quad && t >0,  \label{cont_eq01}\\
   U'_0 (L+t) -G' (t) &= \frac1{c_1}H' ( t) + \frac2{c_1}\sum_{n=1}^{\left[ \tfrac{c_1 t}{2L}\right] } H' \left( t  - \tfrac{2n L}{c_1} \right), \quad && t >0, \label{cont_eq03}
\end{alignat}
\end{subequations}
which by integrating and using \eqref{compatibility} takes the compact form
\begin{equation}\label{final_eq1}
G(t) = \frac{c_1}{c_1 +1} U_0 (L+t)
- \frac2{c_1 +1} \sum_{n=1}^{\left[ \tfrac{c_1 t}{2L}\right] } G \left( t  - \tfrac{2n L}{c_1} \right), \quad t\geq 0,
\end{equation}
under the convention that if $0\leq t<2L/c_1,$ then the sum yields no terms, thus vanishes.


The solution of the above system of equations gives the boundary function $G$ in terms of the known initial condition $U_0.$
Then, the wave is given by \eqref{sol_bou_1} (by replacing $H$ with $G$) for $0\leq x \leq L$, and by \eqref{sol_un_1} for $x>L.$ 

 The latter equation is also used for deriving the backscattered wave at some position $x=L+D,$ with $D>0.$ Thus, by determining the solution $G(t)$ of \eqref{final_eq1}, we will be able to compute the function
\begin{equation}\label{eq_data}
m(t) := u(L+D,t), \quad t>0.
\end{equation}

\begin{theorem}
Let $U_0$ be supported in $(L,\, \infty),$ then the solution of \eqref{final_eq1} is given by
\begin{equation}\label{eq_solution}
\begin{aligned}
G(t) &= \frac{c_1}{c_1 +1} U_0 (L+t) 
- \frac{2c_1}{(c_1 + 1)^2} \sum_{n=0}^{\left[ \tfrac{c_1 t}{2L}-1\right] } \frac{(c_1 - 1)^n}{(c_1 +1)^{n}} U_0 \left(L+ t  - \tfrac{2(n+1) L}{c_1} \right), \quad t \geq 0,
\end{aligned}
\end{equation}
under the convention that if $0\leq t<2L/c_1,$ then the above sum yields no terms, thus vanishes.
\end{theorem}

\begin{proof}
The proof is done by induction. By definition, if $t\in [0, 2L/c_1 ),$ then \eqref{final_eq1} and \eqref{eq_solution} coincide. Let \eqref{eq_solution} be true for $t\in[0,T),$ we will employ \eqref{final_eq1} to show that it is also valid for $t\in [0, T+2L/c_1).$

We consider \eqref{final_eq1} for $t\in [0, T+2L/c_1),$ and as in \eqref{sol_bou_02}, we rewrite it as
\begin{equation}\label{eq_the1}
G(t) = \frac{c_1}{c_1+1} U_0 (L+t)
- \frac2{c_1 +1} \sum_{n=0}^{\left[ \tfrac{c_1 t}{2L}-1\right] } G \left( t  - \tfrac{2(n+1) L}{c_1} \right).
\end{equation}
Since $ t  - \tfrac{2(n+1) L}{c_1} <T+\tfrac{2 L}{c_1}  - \tfrac{2(n+1) L}{c_1} < T+\tfrac{2 L}{c_1}-\tfrac{2 L}{c_1}=T$, the term in the summation can now by computed using \eqref{eq_solution}, thus
\begin{equation}\label{eq_the2}
\begin{aligned}
G \left(t  - \tfrac{2(n+1) L}{c_1}\right) &= \frac{c_1}{c_1 +1} U_0 \left(L+t  - \tfrac{2(n+1) L}{c_1} \right) \\
&\phantom{=}- \frac{2 c_1}{(c_1 +1)^2} \sum_{k=0}^{\left[ \tfrac{c_1 t}{2L} -n-2\right] } \frac{(c_1 - 1)^k}{(c_1 +1)^k} U_0 \left(L+ t  - \tfrac{2(n+k+2) L}{c_1} \right).
\end{aligned}
\end{equation}
Substituting \eqref{eq_the2} into \eqref{eq_the1}, results in
\begin{equation}\label{eq_the3}
\begin{aligned}
G(t) &= \frac{c_1}{1+c_1} U_0 (L+t)
- \frac{2c_1}{(c_1 +1)^2} \sum_{n=0}^{\left[ \tfrac{c_1 t}{2L}-1\right] } U_0 \left(L+t  - \tfrac{2(n+1) L}{c_1} \right) \\
&\phantom{=}+ \frac{4 c_1}{(c_1 +1)^3} \sum_{n=0}^{\left[ \tfrac{c_1 t}{2L}-1\right] } \sum_{k=0}^{\left[ \tfrac{c_1 t}{2L} -n-2\right] } \frac{(c_1 - 1)^k}{(c_1 +1)^k} U_0 \left( L+ t  - \tfrac{2(n+k+2) L}{c_1} \right).
\end{aligned}
\end{equation}

We can simplify the double summation into a single summation, by defining $\lambda = n+k+1,$ and adding ``diagonally" the double sum on $\lambda,$ meaning, by collecting all terms involving $U_0 \left(L+t-\tfrac{2(\lambda+1) L}{c_1}\right).$ In more detail, the double sum in \eqref{eq_the3} is rewritten in the form
\begin{align}\label{eq_sum}
\sum_{n=0}^{\left[ \tfrac{c_1 t}{2L}-1\right] } \sum_{\lambda=n+1}^{\left[ \tfrac{c_1 t}{2L} -1\right] } \left(\frac{c_1 - 1}{c_1 +1}\right)^{\lambda-n-1}  &U_0 \left( L+ t  - \tfrac{2(\lambda+1) L}{c_1} \right) \notag\\
&=\sum_{\lambda=1}^{\left[ \tfrac{c_1 t}{2L}-1\right]} \sum_{n=0}^{\lambda -1}  \left(\frac{c_1 - 1}{c_1 +1}\right)^{\lambda-n-1}  U_0 \left(L+t-\tfrac{2(\lambda +1)L}{c_1}\right)  \notag \\
&=\sum_{\lambda=1}^{\left[ \tfrac{c_1 t}{2L}-1\right]}  \left( \sum_{\mu=0}^{\lambda -1} \left(\frac{c_1 - 1}{c_1 +1}\right)^{\mu} \right)  U_0 \left(L+t-\tfrac{2(\lambda +1)L}{c_1}\right)    \notag \\
&=  -\frac{c_1 +1}{2} \sum_{\lambda=0}^{\left[ \tfrac{c_1 t}{2L}-1\right] }  \left[\left( \frac{c_1 -1}{c_1 +1}\right)^{\lambda} -1 \right] U_0 \left(L+t-\tfrac{2(\lambda+1) L}{c_1}\right). 
\end{align}

Substituting the double sum  in \eqref{eq_the3} by the expression at the RHS of \eqref{eq_sum}, with the change of notation $\lambda = n$,  results in \eqref{eq_solution} for $t\in [0, T+2L/c_1)$, completing the induction.
\end{proof}

\begin{remark}\label{rem-spl-1}
An alternative way for obtaining $G(t)$ follows by the observation that the function $Q(t),$ appearing in the right-hand side of \eqref{sol_bou_3}, can be analysed as 
\begin{subequations}\label{eq_qq}
\begin{alignat}{3}
Q(t) &= \tfrac1{c_1} H'(t), \quad && t\in \left(0, \tfrac{2L}{c_1}\right),\label{eq_qq0}\\
Q(t) &= \tfrac1{c_1} H'(t) + \tfrac2{c_1} H'\left(t-\tfrac{2L}{c_1} \right), \quad && t\in \left(\tfrac{2L}{c_1}, \tfrac{4L}{c_1}\right),\label{eq_qq1}\\
&\phantom{=}\vdots \nonumber\\
Q(t) &= \tfrac1{c_1} H'(t) + \tfrac2{c_1} \sum_{n=1}^N H'\left(t-\tfrac{2nL}{c_1} \right), \quad && t\in \left(\tfrac{2NL}{c_1}, \tfrac{2(N+1)L}{c_1}\right), \label{eq_qqn}
\end{alignat}
\end{subequations}
for $N \in \N_0,$ and with $N=0,$ we mean that no summation is performed.

Thus, we could also compute the solution of \eqref{cont_eq00} by dividing time in appropriate  intervals, as in \eqref{eq_qq}, and compute the boundary functions sectionally. This approach is presented in Appendix \ref{appendix}.
\end{remark}

By substituting \eqref{eq_solution} in \eqref{sol_un_1}  we obtain the measured data \eqref{eq_data} in a d'Alembert form; this reads as the following proposition.
\begin{proposition}
Let $U_0$ be supported in $(L,\, \infty),$ then the measurements at the detector is given by
\begin{equation}\label{data_general1}
m(t) = \frac12 U_0 (L+D+t) + \frac12 U_0 (L+D-t), \quad 0 \leq t \leq D,
\end{equation}
and for $t>D,$ we get
  \begin{equation}\label{data_general2}
\begin{aligned}
m(t) &= \frac12 U_0 (L+D+t) 
  +\frac12\frac{c_1 -1}{c_1 +1} U_0 (L-D+t) \\
  &\phantom{=}-2 c_1 \sum_{n=0}^{\left[ \tfrac{c_1 (t-D)}{2L}-1\right]}\frac{(c_1 -1)^n}{(c_1+1)^{n+2}} U_0 \left(L-D -\frac{2(n+1)L}{c_1}+t\right).
\end{aligned}
\end{equation}
\end{proposition}

\subsubsection{Delta-like short pulse}\label{sec_delta}
In the ideal case of an initial function of the form
\begin{equation}\label{initial_pulse}
U_0 (x) = \left\{\begin{array}{lr}
        1, & \text{for } x=x_0,\\
        0, & \text{for } x \neq x_0,
        \end{array}\right.
\end{equation}
with $x_0 >L,$ the previous result allows us to represent $G,$ (and consequently $m$) at specific time points, as a function of $c_1$ only.

Without loss of generality, we place the source at the detector position, meaning $x_0 = L+D.$ Defining $G^{(k)} := G(D + k\tfrac{2L}{c_1}),\ k\in \N$ and evaluating \eqref{final_eq1}  at  $t=D+k\tfrac{2L}{c_1}, \ k\in \N$, we obtain
\begin{align}
\label{G01-ps}&G^{(0)} = \frac{c_1}{c_1 + 1},  \qquad G^{(1)} =  -\frac{2c_1}{(c_1 + 1)^2}, \\
\label{Gk-ps}&G^{(k)} = - \frac2{c_1 +1} \sum_{n=1}^{\left[ \tfrac{c_1 D}{2L}+k\right] } G^{(k-n)}, \qquad k\ge 2.
\end{align}
Then, the subtraction of the $k$th term from the $(k-1)$th term of \eqref{Gk-ps} yields
\begin{equation}
\label{Gk-ps-it}G^{(k)} = \frac{c_1 -1}{c_1 +1}G^{(k-1)}, \qquad k\ge 2.
\end{equation}
Thus, \eqref{G01-ps} and \eqref{Gk-ps-it} yield
\begin{equation}
G^{(k)} = \begin{cases}\vspace{0.2cm}
        \frac{c_1}{c_1 + 1}, & \text{for } k=0,\\ \vspace{0.2cm}
         -\frac{2c_1 (c_1 -1)^{k-1}}{(c_1 + 1)^{k+1}}, & \text{for } k=1,2,\ldots
        \end{cases}
\end{equation}

Setting $t=0$ in \eqref{data_general1}, we get $m(0)=1,$ as expected. We define
$m^{(k)} := m(2D + k\tfrac{2L}{c_1}),$ and from \eqref{data_general2} we get
\begin{equation}\label{data_delta}
m^{(k)} =  \begin{cases}\vspace{0.2cm}
        \frac12\frac{c_1-1}{c_1 + 1}, & \text{for } k=0,\\ \vspace{0.2cm}
G^{(k)}, & \text{for } k=1,2,\ldots
        \end{cases}
\end{equation}

\subsection{Double-layered medium}

We consider the same scattering problem but now the medium consists of two layers, again supported in $(0,L).$ Let the $j$th layer have wave speed $c_j,$ and length $l_j,$ for $j=1,2,$ such that $l_1 + l_2 = L.$ 

This scattering problem can be decomposed into three sub-problems in the intervals $(0,\, l_2),$ $(l_2,\, L),$ and $(L,+\infty)$ as follows:
 
\begin{description}
\item[IBVP I (half line)]
\begin{subequations}\label{system2_un}
\begin{alignat}{3}
u_{tt} - u_{xx} & = 0,  \quad &&x >L, \, t>0,  \label{eq2_un1}\\
u(x,0) = U_0 (x), \, u_t (x,0) &=0,   \quad &&x >L, \label{ic2_un1}\\ 
u(L,t) = G_0 (t), \, u_x(L,t) &= F_0 (t), \quad &&t>0, \label{bc2_un1}
\end{alignat}
\end{subequations}
where $U_0 : [L,\,+\infty) \rightarrow \R.$
 
\item[IBVP II (finite interval)]
\begin{subequations}\label{system2_bou}
\begin{alignat}{3}
v_{tt} - c_1^2 v_{xx} & = 0,  \quad && l_2 <x <L, \, t>0,  \label{eq2_b1}\\
v(x,0) = v_t (x,0) &=0,   \quad && l_2 <x <L, \label{ic2_b1}\\ 
v(l_2,t) = G_1(t), \, v_x(l_2,t) &= F_1 (t), \quad && t>0, \label{bc2_b1}\\ 
 v(L,t) = H_1 (t), \, v_x(L,t) &= Q_1 (t), \quad &&t>0. \label{bc2_b2}
\end{alignat}
\end{subequations}

\item[IBVP III (finite interval)]
\begin{subequations}\label{system3_bou}
\begin{alignat}{3}
w_{tt} - c_2^2 w_{xx} & = 0,  \quad && 0 <x <l_2, \, t>0,  \label{eq3_b1}\\
w(x,0) = w_t (x,0) &=0,   \quad && 0 <x <l_2, \label{ic3_b1}\\ 
w(0,t) &= 0,  \quad && t>0, \label{bc3_b1}\\ 
 w(l_2,t) = H_2 (t), \, w_x(l_2,t) &= Q_2 (t), \quad &&t>0. \label{bc3_b2}
\end{alignat}
\end{subequations}
\end{description} 
 
The compatibility condition now reads
\begin{equation}\label{compatibility2}
U_0 (L) = G_0 (0) = G_1 (0) = H_1 (0) = H_2 (0)=0.
\end{equation} 
 
The problems \eqref{system2_un}, \eqref{system2_bou}, and  \eqref{system3_bou} are coupled through the boundary -- continuity conditions:

\begin{subequations}\label{cont2_eq0}
\begin{alignat}{2}
G_0 (t) &= u(x,t) \mid_{x=L}= v(x,t) \mid_{x=L}=H_1 (t),  \quad \quad && t>0, \label{cont2_1a}\\
 F_0 (t) &= u_x(x,t) \mid_{x=L} = v_x(x,t) \mid_{x=L}= Q_1 (t) , \quad \quad && t>0, \label{cont2_1b} \\
 G_1 (t) &= v(x,t) \mid_{x=l_2}= w(x,t) \mid_{x=l_2}=H_2 (t),  \quad \quad && t>0, \label{cont2_2a}\\
 F_1 (t) &= v_x(x,t) \mid_{x=l_2} = w_x(x,t) \mid_{x=l_2}= Q_2 (t) , \quad \quad && t>0. \label{cont2_2b}
\end{alignat}
\end{subequations} 
%

The solution of \eqref{system2_un} is given by \eqref{sol_un_1} by replacing $G$ with $G_0.$ Similarly, the solution of \eqref{system3_bou} is given by
\begin{equation}\label{sol3_bou_1}
w(x,t) = \sum_{n=0}^{\left[ \tfrac{c_2 t +x-l_2}{2l_2}\right] } H_2 \left( t +\tfrac{x - (2n+1)l_2}{c_2}\right) - \sum_{n=0}^{\left[ \tfrac{c_2 t -x-l_2}{2l_2}\right] } H_2 \left( t - \tfrac{x + (2n+1)l_2}{c_2}\right),
\end{equation}
using \eqref{sol_bou_1} in this problem. 
 
Using \cite[Eq. (63)]{FokKal22} for $t = \tilde{t}/c_1,$ and $x = \tilde{x} + l_2,$ we get
 \begin{equation}\label{sol2_bou_1}
 \begin{aligned}
v(x,t) &= \sum_{n=0}^{\left[ \tfrac{c_1 t +x-l_2-l_1}{2l_1}\right] } H_1 \left( t +\tfrac{x - l_2 - (2n+1)l_1}{c_1}\right) 
- \sum_{n=0}^{\left[ \tfrac{c_1 t -x + l_2 -l_1}{2l_1}\right] } H_1 \left( t - \tfrac{x -l_2 + (2n+1)l_1}{c_1}\right)  \\
 &\phantom{=}-\sum_{n=0}^{\left[ \tfrac{c_1 t +x-l_2-2l_1}{2l_1}\right] } G_1 \left( t +\tfrac{x - l_2 - 2(n+1)l_1}{c_1}\right) 
 + \sum_{n=0}^{\left[ \tfrac{c_1 t -x + l_2}{2l_1}\right] } G_1 \left( t - \tfrac{x -l_2 +  2n l_1}{c_1}\right),
 \end{aligned}
\end{equation}
 where we have set $h_1 (c_1 t) =H_1 (t),$ and $g_1 (c_1 t) =G_1 (t).$ The initial conditions \eqref{ic2_b1} are clearly satisfied and at the right boundary $x=L = l_1 +l_2,$ we get
  \begin{equation}\label{sol2_bou_2}
 \begin{aligned}
v(L,t) &= \sum_{n=0}^{\left[ \tfrac{c_1 t }{2l_1}\right] } H_1 \left( t -\tfrac{2n l_1}{c_1}\right) - \sum_{n=0}^{\left[ \tfrac{c_1 t}{2l_1}\right]-1 } H_1 \left( t - \tfrac{2(n+1)l_1}{c_1}\right)  \\
 &\phantom{=}-\sum_{n=0}^{\left[ \tfrac{c_1 t -l_1}{2l_1}\right] } G_1 \left( t -\tfrac{ (2n+1)l_1}{c_1}\right) + \sum_{n=0}^{\left[ \tfrac{c_1 t - l_1}{2l_1}\right] } G_1 \left( t - \tfrac{(2n+1) l_1}{c_1}\right) 
 \\
 &= H_1 (t), 
 \end{aligned}
\end{equation}
 and
  \begin{equation}\label{sol2_bou_3}
 \begin{aligned}
v_x(L,t) &= \frac1{c_1}\sum_{n=0}^{\left[ \tfrac{c_1 t }{2l_1}\right] } H'_1 \left( t -\tfrac{2n l_1}{c_1}\right) + \frac1{c_1}\sum_{n=0}^{\left[ \tfrac{c_1 t}{2l_1}\right]-1 } H'_1 \left( t - \tfrac{2(n+1)l_1}{c_1}\right)  \\
 &\phantom{=}-\frac1{c_1}\sum_{n=0}^{\left[ \tfrac{c_1 t -l_1}{2l_1}\right] } G'_1 \left( t -\tfrac{ (2n+1)l_1}{c_1}\right) -\frac1{c_1} \sum_{n=0}^{\left[ \tfrac{c_1 t - l_1}{2l_1}\right] } G'_1 \left( t - \tfrac{(2n+1) l_1}{c_1}\right) 
 \\
 &= \frac1{c_1}H'_1 (t) + \frac2{c_1}\sum_{n=1}^{\left[ \tfrac{c_1 t }{2l_1}\right] } H'_1 \left( t -\tfrac{2n l_1}{c_1}\right) -\frac2{c_1} \sum_{n=0}^{\left[ \tfrac{c_1 t - l_1}{2l_1}\right] } G'_1 \left( t - \tfrac{(2n+1) l_1}{c_1}\right)\\
 &=: Q_1 (t). 
 \end{aligned}
\end{equation}

 At the left boundary $x=l_2,$ we obtain
  \begin{equation}\label{sol2_bou_4}
 \begin{aligned}
v(l_2,t) &= \sum_{n=0}^{\left[ \tfrac{c_1 t -l_1}{2l_1}\right] } H_1 \left( t -\tfrac{(2n+1) l_1}{c_1}\right) - \sum_{n=0}^{\left[ \tfrac{c_1 t -l_1}{2l_1}\right] } H_1 \left( t - \tfrac{(2n+1)l_1}{c_1}\right)  \\
 &\phantom{=}-\sum_{n=0}^{\left[ \tfrac{c_1 t}{2l_1}\right]-1 } G_1 \left( t -\tfrac{ 2(n+1)l_1}{c_1}\right) + \sum_{n=0}^{\left[ \tfrac{c_1 t }{2l_1}\right] } G_1 \left( t - \tfrac{2n l_1}{c_1}\right) 
 \\
 &= G_1 (t), 
 \end{aligned}
\end{equation}
 and
  \begin{equation}\label{sol2_bou_5}
 \begin{aligned}
 \begin{aligned}
v_x (l_2,t) &= \frac1{c_1}\sum_{n=0}^{\left[ \tfrac{c_1 t -l_1}{2l_1}\right] } H'_1 \left( t -\tfrac{(2n+1) l_1}{c_1}\right) +\frac1{c_1} \sum_{n=0}^{\left[ \tfrac{c_1 t -l_1}{2l_1}\right] } H'_1 \left( t - \tfrac{(2n+1)l_1}{c_1}\right)  \\
 &\phantom{=}-\frac1{c_1}\sum_{n=0}^{\left[ \tfrac{c_1 t}{2l_1}\right]-1 } G'_1 \left( t -\tfrac{ 2(n+1)l_1}{c_1}\right) -\frac1{c_1} \sum_{n=0}^{\left[ \tfrac{c_1 t }{2l_1}\right] } G'_1 \left( t - \tfrac{2n l_1}{c_1}\right) 
 \\
 &= -\frac1{c_1} G'_1 (t)-\frac2{c_1}\sum_{n=1}^{\left[ \tfrac{c_1 t }{2l_1}\right] } G'_1 \left( t -\tfrac{2n l_1}{c_1}\right)+  \frac2{c_1}\sum_{n=0}^{\left[ \tfrac{c_1 t -l_1}{2l_1}\right] } H'_1 \left( t -\tfrac{(2n+1) l_1}{c_1}\right) \\
 &=: F_1 (t).
 \end{aligned}
 \end{aligned}
\end{equation}
The function $F_1$ is similar to  $Q_1$ with the roles of $H_1$ and $G_1$ interchanged. 

Finally, we substitute \eqref{sol_un_2} and \eqref{sol_un_3} (for $G_0$ and $F_0$) and \eqref{sol2_bou_2}
-- \eqref{sol2_bou_5} in \eqref{cont2_eq0} to derive 
\begin{subequations}\label{cont2_eq}
\begin{alignat}{3}
  G_0 (t) &= H_1 ( t),  \quad && t>0, \label{cont2_eq1}\\
   U'_0 (L+t) -G'_0 (t) &=  Q_1 (t), \quad && t>0, \label{cont2_eq3}\\
    G_1 (t) &= H_2 (t),  \quad && t>0, \label{cont2_eq4}\\
   F_1 (t) &= Q_2 (t),  \quad && t>0. \label{cont2_eq5}
\end{alignat}
\end{subequations}

We consider the condition \eqref{compatibility2} and we write the above system of equations only for the two unknown functions $G_0$ and $G_1$

\begin{subequations}\label{layer2_final}
\begin{alignat}{3}
G_0 (t) &= \frac{c_1}{c_1 +1} U_0 (L+t) -\frac{2}{c_1 +1} \sum_{n=1}^{\left[ \tfrac{c_1 t }{2l_1}\right] } G_0 \left( t -\tfrac{2n l_1}{c_1}\right)
 + \frac2{c_1 +1} \sum_{n=0}^{\left[ \tfrac{c_1 t - l_1}{2l_1}\right] } G_1 \left( t - \tfrac{(2n+1) l_1}{c_1}\right),\label{layer2_final1}
\\
G_1 (t) &= \frac{2c_2}{c_2 +c_1} \sum_{n=0}^{\left[ \tfrac{c_1 t - l_1}{2l_1}\right] } G_0 \left( t - \tfrac{(2n+1) l_1}{c_1}\right) -\frac{2c_2}{c_2 +c_1} \sum_{n=1}^{\left[ \tfrac{c_1 t }{2l_1}\right] } G_1 \left( t -\tfrac{2n l_1}{c_1}\right) 
-\frac{2c_1}{c_2 +c_1} \sum_{n=1}^{\left[ \tfrac{c_2 t }{2l_2}\right] } G_1 \left( t -\tfrac{2n l_2}{c_2}\right). \label{layer2_final2}
\end{alignat}
\end{subequations}

The system of equations \eqref{layer2_final} is the analogue to \eqref{final_eq1} for the double-layered medium. However, the two unknown functions are not only coupled but also evaluated at different time intervals.

\begin{remark}\label{rem-spl-2}
 In order to derive a formula for the measurements, it is sufficient to obtain $G_0$. This is doable using the methodology described in Remark \ref{rem-spl-1} to solve \eqref{layer2_final}, which yields the solution in the form
 \begin{equation}
 G_0(t)= \sum_{k_1=0}^{\left[ \tfrac{c_1 t }{2l_1}\right] }\sum_{k_2=0}^{\left[ \tfrac{c_2 t }{2l_2}\right] } A^{(2k_1,2k_2)} U_0 \left( L+t  - k_1\tfrac{2l_1}{c_1}- k_2\tfrac{2l_2}{c_2}\right),
 \end{equation}
where $A^{(2k_1,2k_2)}$ are real constants depending on $c_1$ and $c_2$.

In this work we omit the elaborate analysis that should be followed in order to determine all the coefficients $A^{(2k_1,2k_2)}$ in the different physical setups. Instead, we present
\begin{align*}
A^{(0,0)} &= \frac{c_1}{1+c_1} \qquad A^{(2,0)} =\frac{2c_1(c_2 - c_1)}{(1+c_1)^2(c_1+c_2)} \\
 A^{(2,2)}& =-\frac{8c_1^2 c_2 }{(1+c_1)^2(c_1+c_2)^2},
\end{align*}
which correspond to the information obtained by the three major reflections, at $x = L, \, l_2,$ and $0,$  respectively.

In fact, for the inverse problem, only the first two coefficients are needed, since these two constraints are enough for obtaining the two unknowns $c_1$ and $c_2$. The third coefficient serves as an assurance that the boundary at $x=0$ is indeed totally reflecting the wave, namely that the Dirichlet condition vanishes at $x=0$. 
\end{remark}

In what follows we present the treatment of the above problem for the case of a delta-like initial function as in Section \ref{sec_delta} in order to simplify the upcoming calculations. 

Let $U_0$ be of the form \eqref{initial_pulse} with $x_0 = L+D.$ Following the above remark, we are interested in computing $G_0$ at the time steps $t = D, \, D+ \tfrac{2l_1}{c_1},$  corresponding to the reflections by the boundaries at the positions $x = L,$ and $l_2,$ respectively. Multiple reflections can also be computed but are not necessary for the inverse problem. We define $G_j^{(k_1 ,k_2)} := G_j(D + k_1\tfrac{l_1}{c_1}+ k_2\tfrac{l_2}{c_2}),$ for $j=0,1$ and $k_j \in \Z.$

We evaluate \eqref{layer2_final1}  at the time steps $t=D$ and $t=D+\tfrac{2l_1}{c_1},$ meaning for $(k_1, k_2) = (0,0)$ and $(k_1, k_2) = (2,0),$ respectively, and we subtract the two formulas. Then, we obtain
\begin{equation}\label{eq_2l_dif1}
G_0^{(2,0)}  - G_0^{(0,0)} = -\frac{c_1}{c_1+1} - \frac{2}{c_1+1} G_0^{(0,0)} +\frac{2}{c_1+1} G_1^{(1,0)}.
\end{equation}

Given that $G_0$ and $G_1$ refer to boundary functions, we recall from \eqref{system2_un} and \eqref{system2_bou} that
\begin{equation}\label{conditions}
\begin{aligned}
G_0 (t) &= 0, \quad \text{for} \quad t<D, \\
G_1 (t) &= 0, \quad \text{for} \quad t<D+\tfrac{l_1}{c_1},
\end{aligned}
\end{equation}
meaning
\begin{equation}\label{conditions_test}
\begin{aligned}
G_0^{(k_1,k_2)} &= 0, \quad \text{for} \quad k_1 <0 \text{ and } k_2 < 0 , \\
G_1^{(k_1,k_2)} &= 0, \quad \text{for} \quad k_1 < 1 \text{ and }  k_2 < 0.
\end{aligned}
\end{equation}

We set $t=D$ in \eqref{layer2_final1}--\eqref{layer2_final2} and using \eqref{initial_pulse} and \eqref{conditions} we get
\begin{equation}\label{eq_2l_g0}
G_0^{(0,0)} = \frac{c_1}{c_1 +1}, \quad \text{and} \quad G_1^{(0,0)} = 0.
\end{equation}
In addition, setting $t=D+\tfrac{l_1}{c_1}$ in \eqref{layer2_final2}, we see that only the first sum contributes a non-zero term and precisely
\begin{equation}\label{eq_2l_g1}
G_1^{(1,0)} = \frac{2 c_2}{c_2 +c_1} G_0^{(0,0)}.
\end{equation}

Then, employing \eqref{eq_2l_g0} and \eqref{eq_2l_g1} in \eqref{eq_2l_dif1}, we obtain 
\begin{equation}\label{eq_g0_l1_final}
G_0^{(2,0)}  = \frac{2 c_1 (c_2 - c_1)}{(c_1 +1)^2 (c_2 +c_1)}. 
\end{equation}

The above procedure can be continued in order to compute the terms describing multiple reflections. For example, starting from
\begin{equation}\label{eq_2l_dif1b}
G_0^{(4,0)}  - G_0^{(2,0)} =  - \frac{2}{c_1+1} G_0^{(2,0)} +\frac{2}{c_1+1} G_1^{(3,0)},
\end{equation}
after some straightforward calculations we compute
\begin{equation}
G_0^{(4,0)} = -\frac{(c_1-1)(c_2-c_1)}{(c_1+1)(c_2+c_1)} G_0^{(2,0)}.
\end{equation}
In general, we get
\begin{equation}
G_0^{(2k+2,0)} = -\frac{(c_1-1)(c_2-c_1)}{(c_1+1)(c_2+c_1)} G_0^{(2k,0)}, \qquad  k=1,2,\ldots,
\end{equation}
which, together with \eqref{eq_g0_l1_final}, yields
\begin{equation}
G_0^{(2k,0)} = 2c_1 \frac{(-1)^{k+1} (c_1-1)^{k-1}(c_2-c_1)^k}{(c_1 +1)^{k+1}(c_2 +c_1)^k}, \quad k=1,2,\ldots
\end{equation}

\begin{remark}
Recalling the discussion in Remark \ref{rem-spl-2}, we note that the major peak induced by the reflection at $x=0$ is given by evaluating $G_0$ at $ t=D+ \tfrac{2l_1}{c_1}+ \tfrac{2l_2}{c_2},$ namely
\begin{equation}\label{eq_g0_l1_l2}
G_0^{(2,2)} = - \frac{8c_1^2 c_2}{(c_2 +c_1)^2 (c_1 + 1)^2}.  
\end{equation} 
The derivation of this formula is given in Appendix \ref{appendix-b}.
\end{remark}

The terms in \eqref{eq_2l_g0} (first equation), \eqref{eq_g0_l1_final} and \eqref{eq_g0_l1_l2} summarize the single/major reflections from the boundary interfaces for delta-like pulse of the form \eqref{initial_pulse} and describe the major terms in the representation of the back-reflected wave.

We have derived all the ingredients needed for the following proposition considering the measurements at $x=L+D$:
\begin{proposition}
Let $m^{(n)} := m\left(2D + \sum_{j=1}^n\tfrac{2 \ell_j}{c_j}\right).$ Then, the three major peaks are given by
\begin{equation}\label{data_delta2}
m^{(n)} =\begin{cases} \vspace{3mm}
        \frac12\frac{c_1-1}{c_1 + 1}, & \text{for } n=0,\\ \vspace{2mm}
G_0^{(2,0)}, & \text{for } n=1,\\ \vspace{0.2cm}
G_0^{(2,2)}, & \text{for } n=2.
        \end{cases}
\end{equation}
where $G_0^{(2,0)}$ and $G_0^{(2,2)}$ are given by \eqref{eq_g0_l1_final} and \eqref{eq_g0_l1_l2}, respectively.
\end{proposition}
\begin{proof}
The proof is straightforward by employing \eqref{sol_un_1} in \eqref{eq_data}, and evaluating the resulting expression at $t=2D + \sum_{j=1}^n\tfrac{2 \ell_j}{c_j}, \ \ n=0,1,2$.
\end{proof}

\section{The inverse scattering problem}\label{sec_inverse}

In this section we examine the numerical solution of the inverse problem, to recover the properties of the medium from the knowledge of the initial function and the measured data. We consider the setup as described in Section \ref{sec_delta} for a delta-like $U_0$ and back-reflected data either full (as defined in \eqref{data_delta}) or phaseless, meaning $m(t) = |u(L+D,t)|.$

\subsection{Single-layered medium}

The data consists of peaks with varying heights at specific time steps. We ignore the one corresponding to the non-reflected initial wave (see for example the most left peak in \autoref{fig3}). Let us denote by $h_k,$ the height of the $k-$th peak appearing at time $t_k,$ for $k=1,2,\ldots .$ We have more than enough information to reconstruct just two unknowns: the wave speed $0<c_1<1$ and the length $L>0.$ 

The peaks appear at the following time steps
\begin{equation}
\begin{aligned}
t_1 &= 2D,\\
t_2 &= 2D + \frac{2L}{c_1},\\
&\vdots \\
t_k &= 2D + (k-1)\frac{2L}{c_1}, \quad k = 1,2,\ldots
\end{aligned}
\end{equation}
Thus, the difference of the first two equations gives us the ratio
\begin{equation}
\frac{L}{c_1} = \frac{t_2 - t_1}{2}. 
\end{equation}
The wave speed will be recovered from the equations of the heights and it will be substituted in the above equation, in order to recover also the length. In what follows, we consider both the cases of full and phaseless data.

\subsubsection{Full-data}

The solution of the inverse problem in this case is trivial. The information of the first peak, namely 
$m^{(0)} = h_{1}$, is sufficient. Hence, we obtain 
the unique solution:
\begin{equation}\label{c1-sol-1}
c_1 = \frac{1+2 h_1}{1-2 h_1}.
\end{equation}

\subsubsection{Phaseless-data}

In this case, we have to solve $|m^{(0)}| = |h_1 |$ which admits two solutions, in general.
Recalling that $c_1 <1,$ then the peak corresponding to the reflected wave from the first boundary will have negative sign, namely $h_1 <0$ (see the second peak in \autoref{fig3}). Thus, in view of \eqref{data_delta} the equation to be solved takes the form:
$$\frac{1}{2}\frac{c_1-1}{c_1+1}=-|h_1|,$$
which yields the unique solution
$$c_1=\frac{1-2 |h_1|}{1+2 |h_1|},$$
which is identical to \eqref{c1-sol-1}.

As we will see later, excluding one solution even for the double-layered medium is not possible and additional information is needed.

\subsection{Double-layered medium}

We are interested in reconstructing the wave speeds $0<c_j<1$ and the lengths $\ell_j >0,$ for $j=1,2.$ 

The major peaks appear at the following time steps
\begin{equation}\label{eq_time}
\begin{aligned}
t_1 &= 2D,\\
t_2 &= 2D + \frac{2 \ell_1}{c_1},\\
t_3 &= 2D + \frac{2 \ell_1}{c_1}+ \frac{2 \ell_2}{c_2}.
\end{aligned}
\end{equation}
As in the single-layered case, the combination of the above equations gives
\begin{equation}\label{eq_ratio}
\frac{\ell_j}{c_j} = \frac{t_{j+1} - t_j}{2}, \quad j=1,2. 
\end{equation}

\subsubsection{Full-data}

Equating the amplitudes, we obtain
\begin{subequations}
\begin{alignat}{3}
\frac{c_1 -1}{2(c_1 +1)} &= h_1, \\
\frac{2 c_1 (c_2 - c_1)}{(c_1 +1)^2 (c_2 +c_1)} &= h_2, 
\end{alignat}
\end{subequations}
which admits the unique solution
\begin{subequations}
\begin{alignat}{3}
\label{c1-sol} c_1&= \frac{1+2 h_1}{1-2 h_1},\\
\label{c2-sol} c_2&= \frac{(1+2 h_1) \left(4 h_1^2-2 h_2-1\right)}{(1-2 h_1) \left(4 h_1^2+2 h_2-1\right)}.
\end{alignat}
\end{subequations}

\subsubsection{Phaseless-data}

We observe that $h_1 <0$ and $h_3 <0,$ but the sign of $h_2$ is determined by the term $c_2 -c_1$ which is unknown. Hence, $c_1$ is uniquely recovered by \eqref{c1-sol}, and $c_2$ is satisfying
\begin{equation}
\frac{2 c_1 ( c_2 - c_1) }{(c_1 +1)^2 (c_2 +c_1)} =  \pm |h_2|, \label{eq_less2a}
\end{equation}
 which admits two solutions, depending on the sign of $h_2,$  namely
 \begin{equation}
 \label{c2-sol-pm} c_2= \frac{(1+2 h_1) \left(4 h_1^2\mp 2 |h_2|-1\right)}{(1-2 h_1) \left(4 h_1^2\pm 2 |h_2|-1\right)}.
 \end{equation}
 
 \begin{remark}
 One idea to eliminate one of the above two solutions is to check them against the third measurement, $m^{(2)}$, namely to substitute \eqref{c1-sol} and \eqref{c2-sol-pm} in
 \begin{equation}
 - \frac{8c_1^2 c_2}{(c_2 +c_1)^2 (c_1 + 1)^2} = h_3, \label{eq_less2b}
\end{equation}  
producing an additional constraint in order to determine the sign of $h_2$. However, this procedure results in
\begin{equation}\label{h3-h1h2}
h_3  =   \frac{2 h_2^2}{1-4 h_1^2} -\frac{1-4 h_1^2}{2},
\end{equation}
which yields no extra information, since the RHS of \eqref{h3-h1h2} is independent of the sign of $h_2$.
%
%
%
 \end{remark}
 
The extra information comes from the nature of the scattering problem. Recall that
\begin{equation}\label{eq_length}
\sum_{j=1}^2 \ell_j = L,
\end{equation}
where the right-hand side can be found from the first equation of \eqref{eq_time}, since both source and detector positions (same here) are known, meaning that the sum $L+D = L + t_1/2$ is given. The last equation, considering  \eqref{eq_ratio}, leads to the  equation
\begin{equation}\label{eq_L}
\sum_{j=1}^2 (t_{j+1} - t_j) c_j = 2L.
\end{equation}
Thus, the correct $c_2$ is the one that solves \eqref{eq_less2a} and satisfies \eqref{eq_L}. Of course, in the double-layer case, one can solve directly \eqref{eq_L} since $c_1$ is already reconstructed. However, this is not true for more layers as we are going to see later.

\subsection{Multi-layered medium}\label{sec_multi}

As in the previous cases, we aim to  reconstruct the wave speeds $0<c_j<1$ and the lengths $\ell_j >0,$ for $j=1,\ldots,N$ of a $N-$layered medium. The data consists of $N+1$ peaks with amplitudes $h_j$ at points $t_j,$ for $j=1,\ldots,N+1.$

The equations \eqref{eq_ratio}, \eqref{eq_length} and \eqref{eq_L} still hold by simply replacing the index $2$ with $N$. 

The system of equations to be solved for the wave speeds admit the general form
\begin{equation*}
\begin{aligned}
\rho_1 ( c_1) :=\frac{c_1 -1}{2(c_1 +1)} &= h_1, \\
\rho_2 (c_1, c_2) := \frac{2 c_1 (c_2 - c_1)}{(c_1 +1)^2 (c_2 +c_1)} &= h_2, \\
\rho_3 (c_1,c_2,c_3) &= h_3, \\
&\vdots \\
\rho_N (c_1,\ldots,c_N) &= h_N, 
\end{aligned}
\end{equation*}
The form of the functions $\rho_j,$ for $j=3,\ldots,N$ can be found by considering the corresponding direct problem for a $N-$layered medium. We omit the explicit formulas above for the sake of presentation.  As discussed in the previous section, the last measurement yields the equation $\rho_{N+1} (c_1,\ldots,c_N) = h_{N+1}$ which does not provide any additional information, thus it is neglected in the formation of the above system of equations.

\subsubsection{Full-data}

Following the same procedure as in the double-layered medium, we solve the first equation of the above system, obtaining a unique value for $c_1$. Then, we substitute $c_1$ in the second equation and solve for $c_2,$ obtaining a unique solution. This simple procedure yields the unique values of the wave speeds $\big\{c_j\big\}_{j=1}^N$.

 \subsubsection{Phaseless-data}

In this case, we know that $h_1<0$ and $h_{N+1}<0.$ Thus, we have to solve
  \begin{equation}\label{final_system}
\begin{aligned}
\rho_1 (c_1) &= h_1,\\
|\rho_2 (c_1, c_2)|  &= |h_2|, \\
|\rho_3 (c_1,c_2,c_3)| &= | h_3|, \\
&\vdots \\
|\rho_N (c_1,\ldots,c_N)| &= | h_N |.
\end{aligned}
\end{equation}

\begin{figure}
\begin{center}
\includegraphics[width=0.9\textwidth]{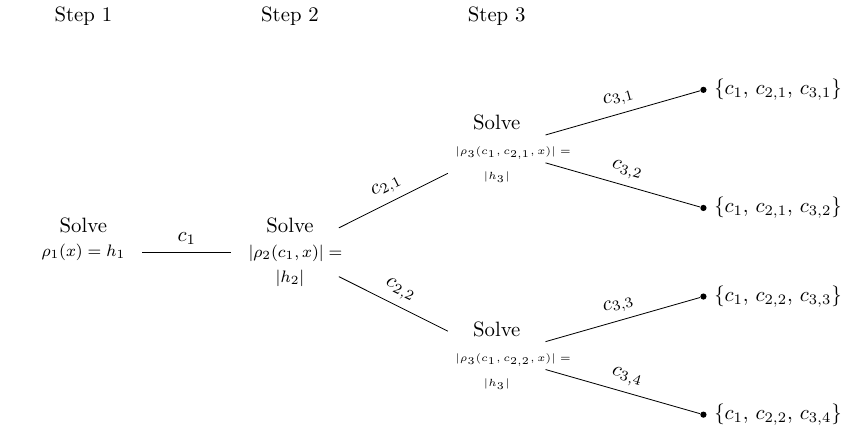}
\caption{The iterative scheme for a $3-$layered medium with length $L$. If all solutions are in the $(0,1)-$interval, we obtain 4 sequences of possible wave speeds. In Step 4, the output of the algorithm is the sequence whose elements satisfy $\sum_{j=1}^3 (t_{j+1}-t_j) c_j = 2L$.}\label{fig_algorithm}
\end{center}
\end{figure}

We propose the following iterative scheme:

\begin{description}
\item{\textbf{Step $\bm 1$:}} Solve $\rho_1 (c_1)  = h_1,$ for $c_1.$ 
\item{\textbf{Step $\bm 2$:}} Substitute $c_1$ in  $|\rho_2 (c_1, c_2)|  = |h_2|,$ and solve it for $c_2,$ obtaining two solutions $c_{2,1}$ and $c_{2,2}.$
\item{\textbf{Step $\bm 3$:}} For the values $c_{2,j_2}\in (0,1), \ \ j_2=1,2,$ solve the equations
\[
|\rho_3 (c_1,c_{2,j_2},c_3)| = | h_3|, 
\]
to obtain $c_3,$ namely maximum four solutions $c_{3,j_3},$ for $j_3=1,2,3,4.$
\item{\quad\vdots} 
\item{\textbf{Step $\bm{N}$:}} Substitute the values $c_{N-1,j_{N-1}}\in (0,1),\ \  j_{N-1}=1,2,\ldots, 2^{N-2}$,  in the last equation of \eqref{final_system} and solve the (at most)  $2^{N-2}$ equations
\[
|\rho_N(c_1,c_{2,j_2},c_{3,j_3},\ldots, c_N)| = | h_N|,
\]
for $c_N,$ to obtain maximum $2^{N-1}$ possible solutions.
\item{\textbf{Step $\bm{N+1}$:}} The set of possible solutions consists of the sequences $\big\{ c_1, c_{2,j_2},\ldots,  c_{N,j_N}\big\},$ with $j_n=1,2,\ldots,2^{n-1}$. The solution  of the inverse problem is the sequence that satisfies \eqref{eq_L}; see \autoref{fig_algorithm} for an illustration.
\end{description}

\section{Numerical examples}\label{sec_numerics}

In the first part, we implement the presented formulas of the solution of the direct problems for single- and double-layered media. Then, we solve the inverse problem to reconstruct the material parameters of a $4-$layered medium using the iterative scheme presented in Section \ref{sec_multi}.

\subsection{The direct scattering problem}

We consider the formulas of Section \ref{sec_direct} for different coefficients and initial functions $U_0.$ In addition, the numerical solutions are compared with the ones obtained from the classical finite difference scheme for piecewise constant wave speed. 

\subsubsection{Single-layered medium}

The solution $W(x,t)$ of the direct problem is given in the following analytical form which provides an extention of the d'Alembert formula for the single-layered medium:
\begin{equation}\label{sol_r}
W(x,t) =  \begin{cases}
  v(x,t), & 0<x \leq L, \, t>0, \\
  u(x,t), & x> L, \, t>0,
\end{cases}
\end{equation}
where $v(x,t)$ is given by \eqref{sol_bou_1} (replace $H$ by $G$) for $0<x \leq L$ and $u(x,t)$ is given by \eqref{sol_un_1} for $x>L.$ The function $G$ is given by \eqref{eq_solution}. 

\begin{figure}
    \centering
    \begin{subfigure}[t]{0.5\textwidth}
        \centering
        \includegraphics[height=2in]{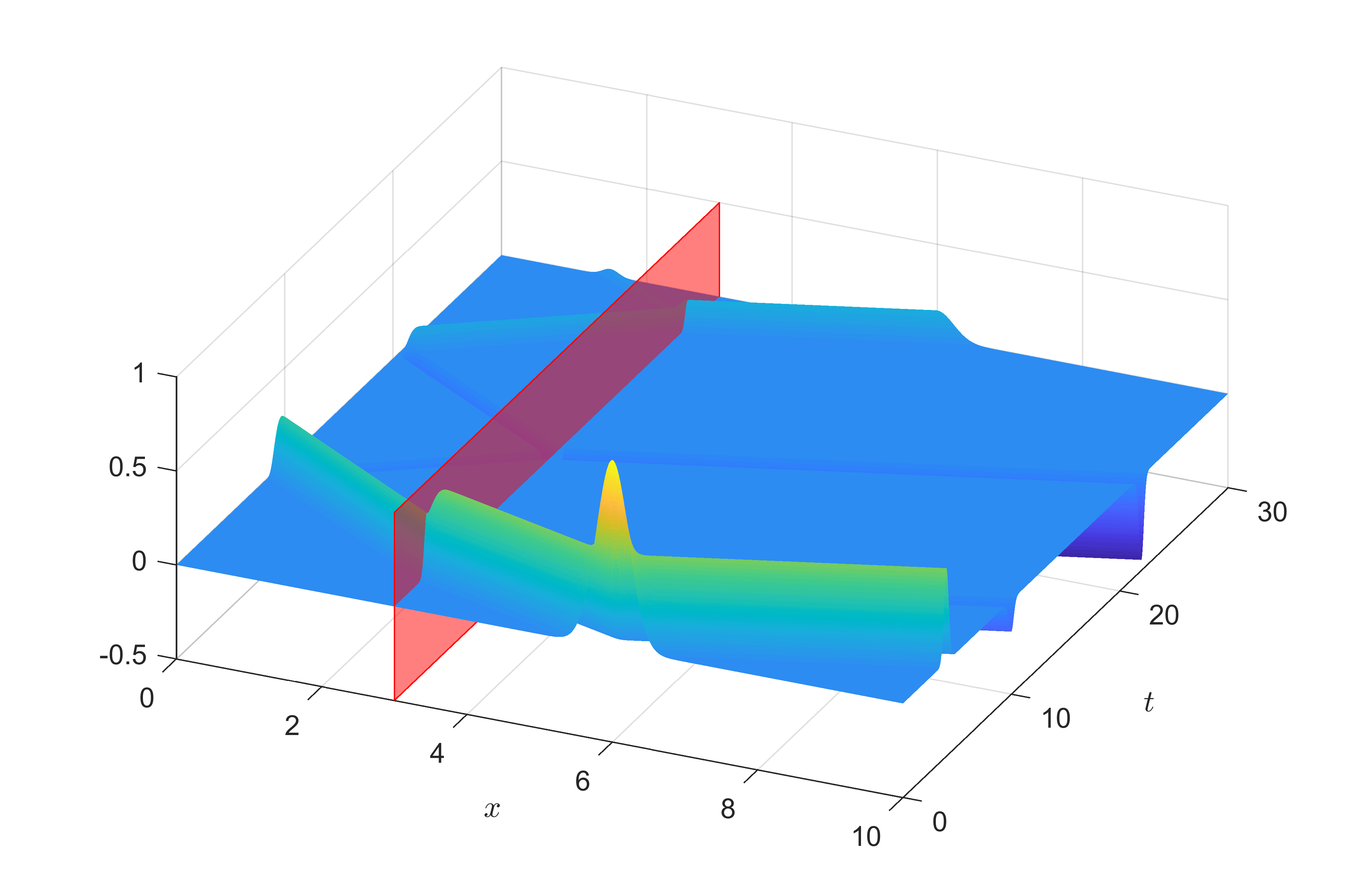}
    \end{subfigure}%
    ~ 
    \begin{subfigure}[t]{0.5\textwidth}
        \centering
        \includegraphics[height=2in]{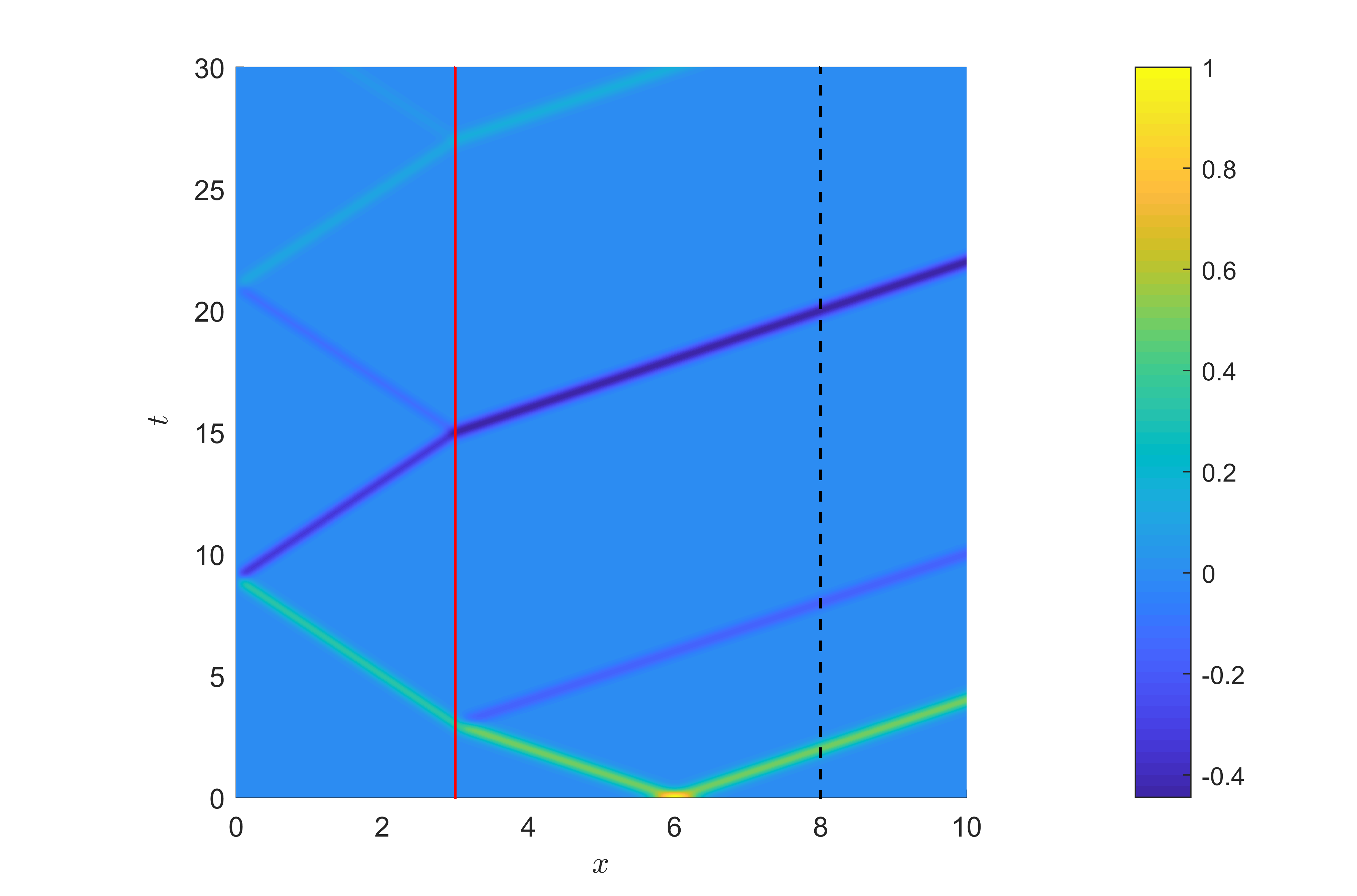}
    \end{subfigure}
    \caption{The function $W$ given by \eqref{sol_r} for the first example (left) and its projection in the $x-t$ plane (right). The plane (in red) identifies the boundary at $x=3.$}\label{fig1}
\end{figure}

In the first example, we set $L=3,$ and $c_1 = 1/2.$ We consider an initial function of the form
\begin{equation}\label{eq_initial}
U_0 (x) = e^{-10(x-x_0)^2}, 
\end{equation}
for $x_0=6,$ and we plot the function $W$ for $x\in [0,10]$ and $t \in [0,30],$ in \autoref{fig1}.

In the second example, we set $L=2,$ and $c_1 = 3/8.$ We consider the continuous but non-differentiable initial condition
\[
U_0 (x) =  \begin{cases}
  0, & x\leq 3\pi/2, \\
  \tfrac{\cos(x)}{x-3}, & x> 3\pi/2.
\end{cases}
\]

\begin{figure}
    \centering
    \begin{subfigure}[t]{0.5\textwidth}
        \centering
        \includegraphics[height=2in]{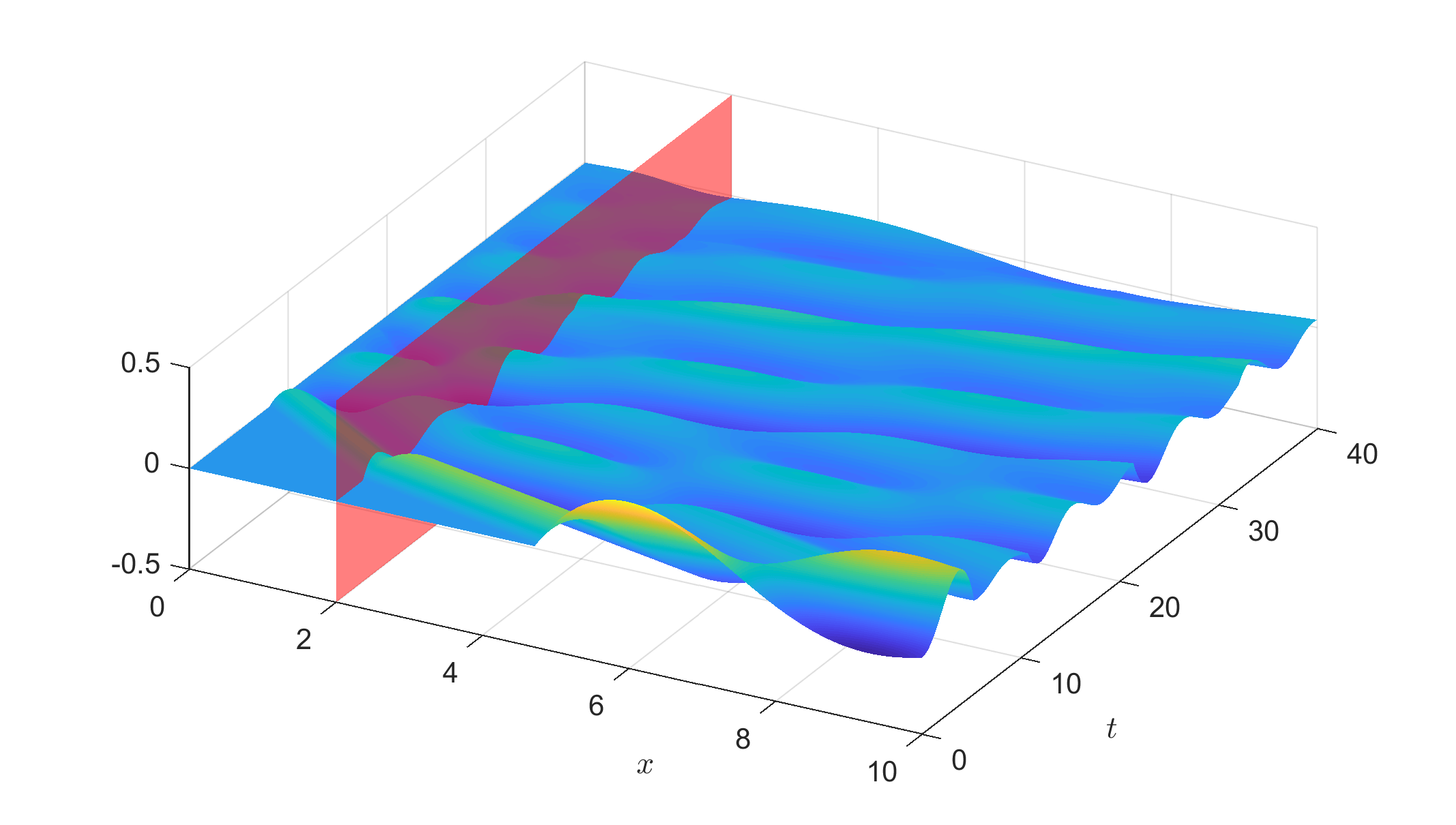}
    \end{subfigure}%
    ~ 
    \begin{subfigure}[t]{0.5\textwidth}
        \centering
        \includegraphics[height=2in]{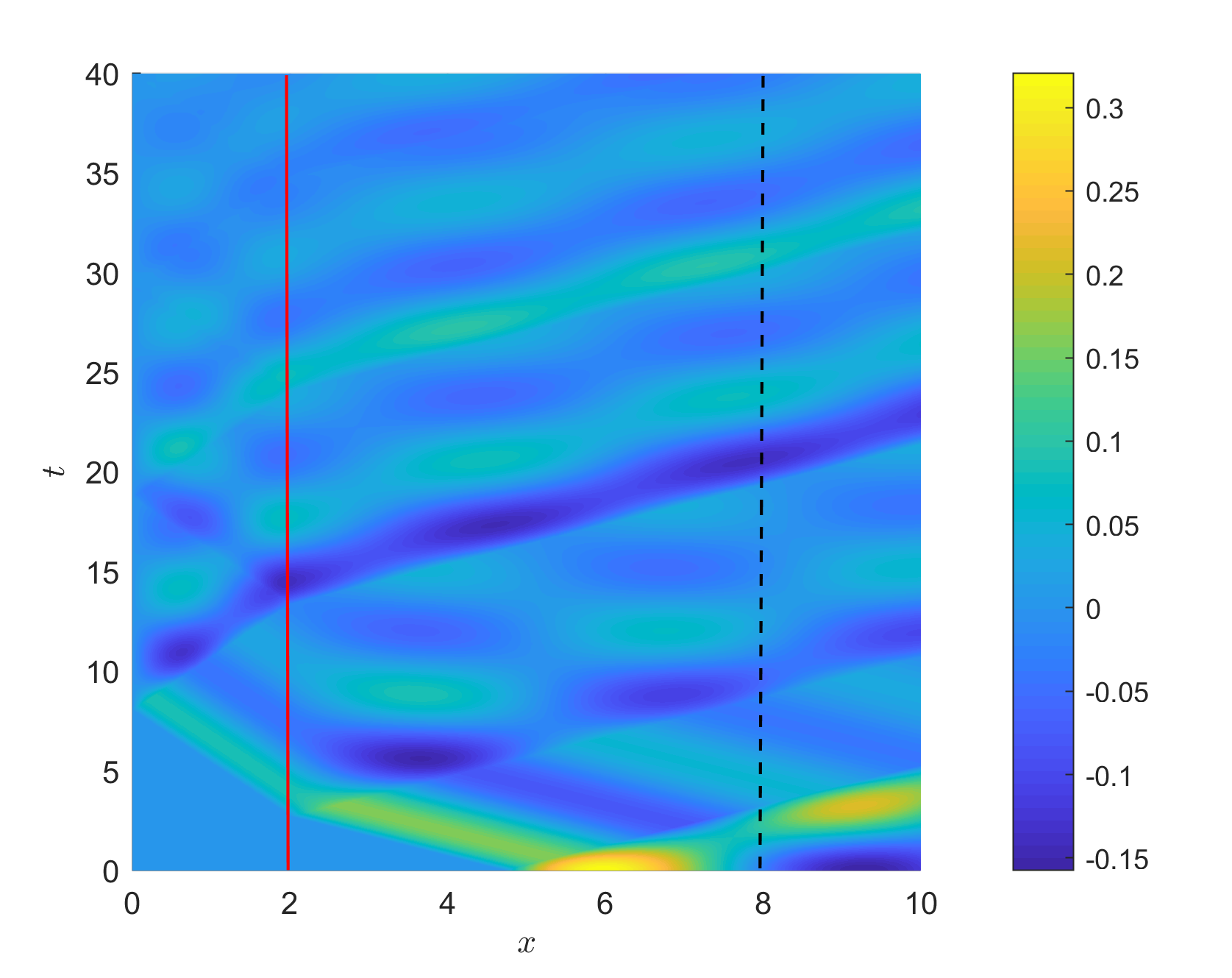}
    \end{subfigure}
    \caption{The function $W$ given by \eqref{sol_r} for the second example (left) and its projection in the $x-t$ plane (right). The plane (in red) identifies the boundary at $x=2.$}\label{fig2}
\end{figure}

The function $W,$ for $x\in [0,10]$ and $t \in [0,40],$ is plotted in  \autoref{fig2}.

\begin{figure}
    \centering
    \begin{subfigure}[t]{0.5\textwidth}
        \centering
        \includegraphics[height=2in]{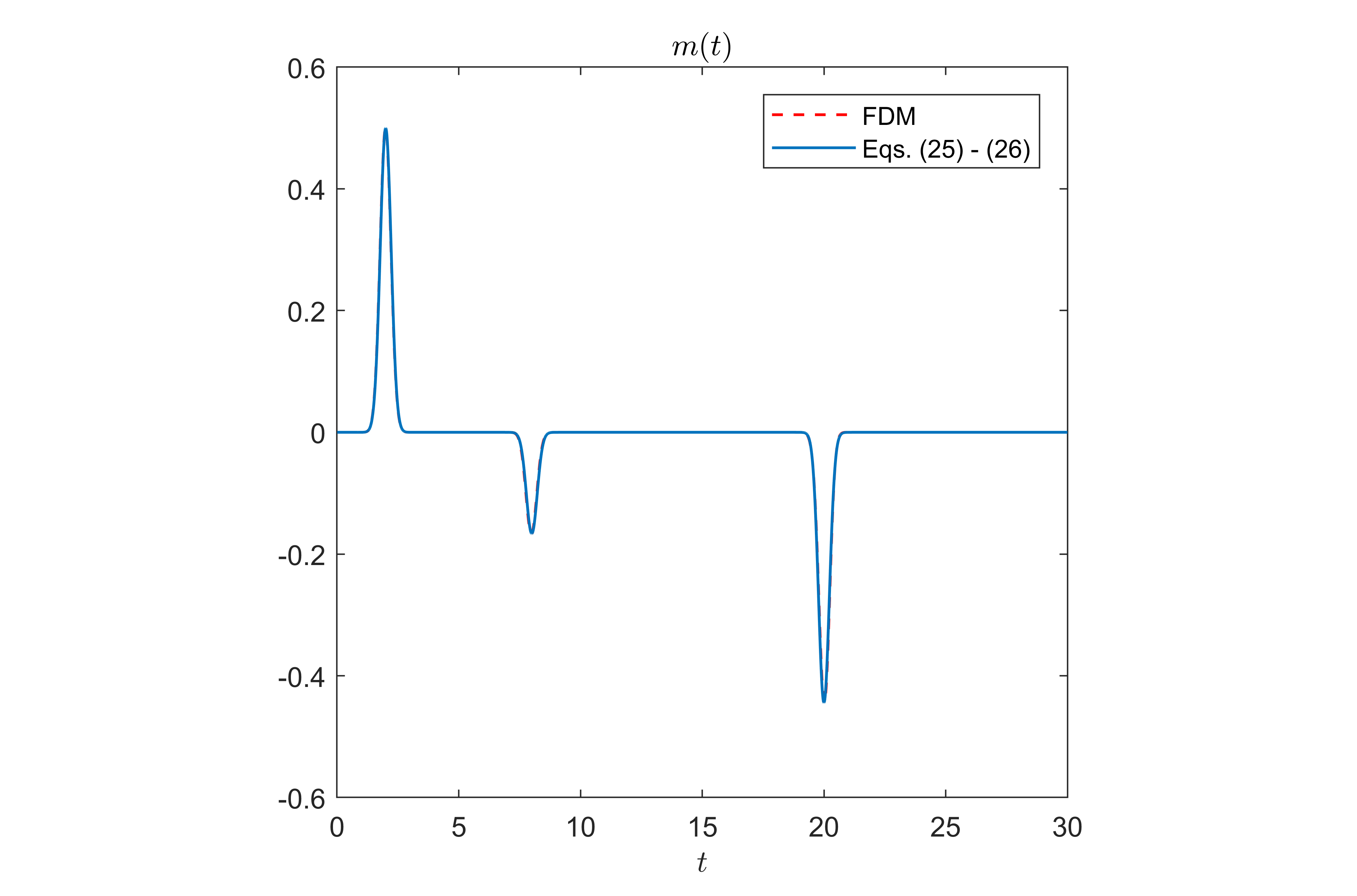}
    \end{subfigure}%
    ~ 
    \begin{subfigure}[t]{0.5\textwidth}
        \centering
        \includegraphics[height=2in]{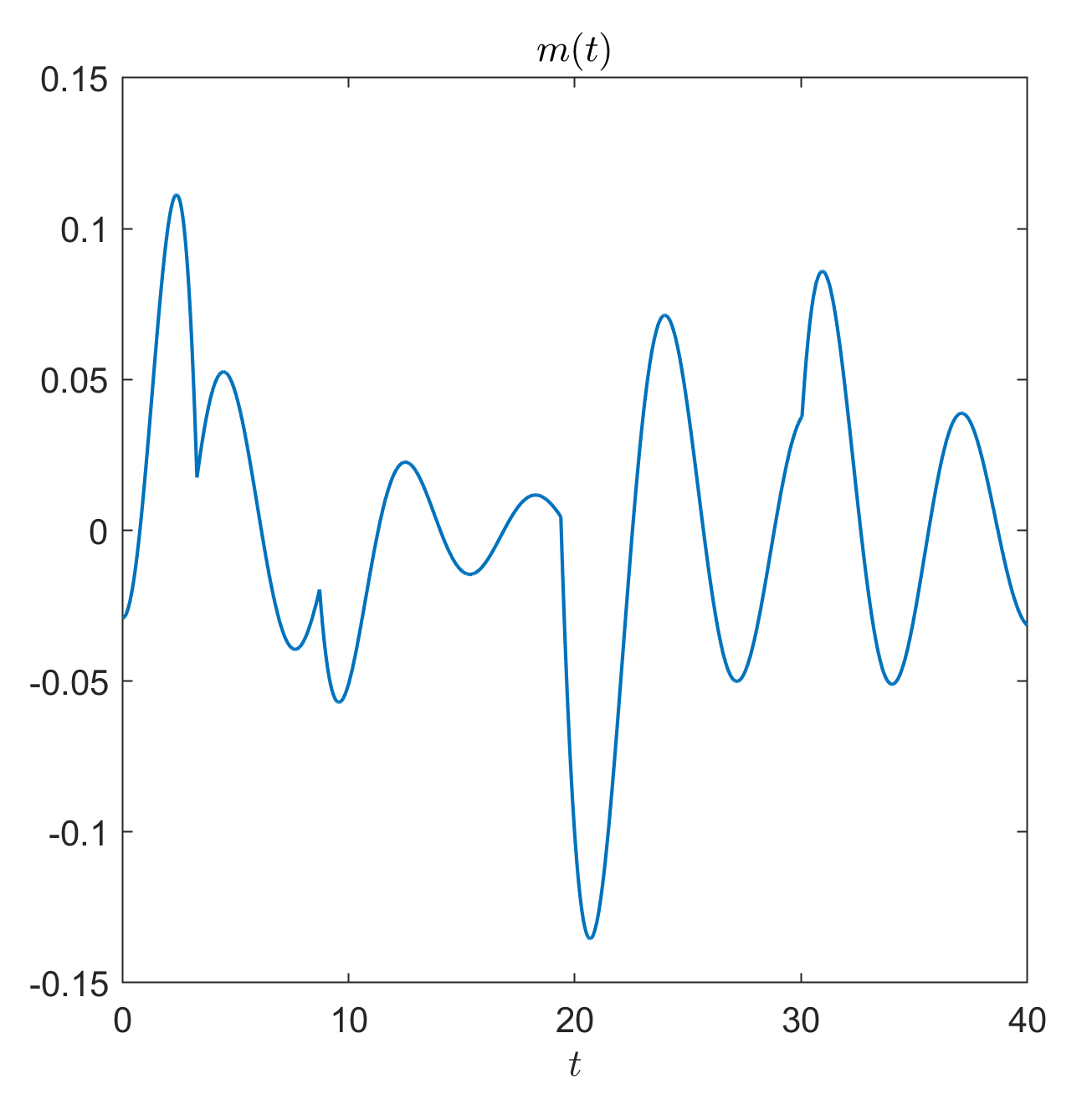}
    \end{subfigure}
\caption{The cross-section of $u$ at $x=8,$ corresponding to the data  $m$ (solid blue line), for the first example (left) and for the second example (right). For the first example this analytical solution is also compared with the numerical solution using FDM (red dotted line).}\label{fig3}
\end{figure}

In \autoref{fig3}, we present the cross-section of $u$ at $x=L+D,$  meaning the data function $m(t),$ for the first example (left) and the second example (right). The position $x=L+D=8,$ is identified with a black dotted line in \autoref{fig1} and \autoref{fig2}.
In addition, we plot the numerical solution of the initial boundary value problem using the finite difference method (FDM) for piecewise constant wave speed. The two solutions are matched perfectly. In the second example, the lack of differentiability of the initial condition prevent us from using the typical FDM. This is another advantage of the proposed scheme since no special discretization is needed.

It is helpful for the corresponding inverse problem (to be addressed later) to discuss the form of the measured data. In the left picture of  \autoref{fig3}, the first peak located at $t=2$ refers to the right-going wave ($x\rightarrow \infty$) originated from $x=6$ and measured at $x=8,$ at the exterior region with $c_0 =1.$ The second and the third peak are the ones of interest, showing the back-reflected waves from the first (at $x=3$) and the second boundary (at $x=0$) of the medium, respectively.

\subsubsection{Double-layered medium}

\begin{figure}
    \centering
        \includegraphics[width=0.8\textwidth]{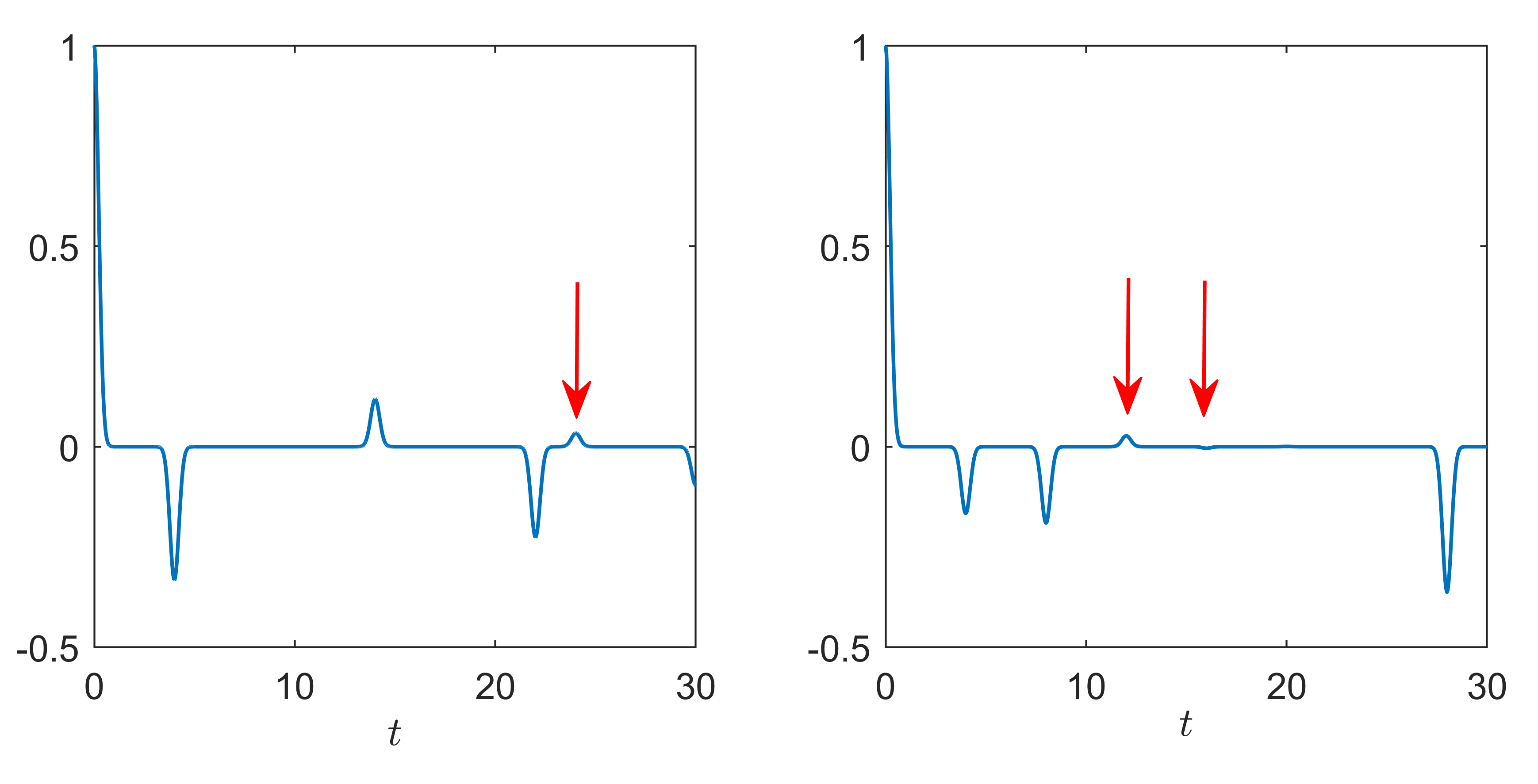}
    \caption{The function $u(5,t),$ for $t\in[0,30]$ of the third example. In the left picture 
    the minor peaks appear after the major peaks, whereas in the right the minor peaks appear  in between the major peaks; the red arrows point at the minor peaks.
    }\label{fig_2layer}
\end{figure}

In the third example, we consider a double-layered medium with $\ell_1 = 1$ and $\ell_2 = 2,$ resulting in $L=3.$ The initial function is given by \eqref{eq_initial} for $x_0 = 5.$ We present $u(5,t)$ for $t\in [0,30]$ for two different sets of wave speeds in \autoref{fig_2layer}. We choose once $(c_1,\,c_2) = (1/2,\,1/5)$ and then $(c_1,\,c_2) = (1/5,\,1/2).$ In the first case, the condition $1<\frac{c_1 l_2}{c_2 l_1} <2,$ see Appendix \ref{appendix-b}, is satisfied. We mark with a red arrow the peaks corresponding to multiple reflected waves (minor peaks). In the first case, all major peaks appear before the minor ones whereas in the second case we observe minor peaks before the last reflection from the boundary at $x=0.$

\subsection{The inverse scattering problem}

In the fourth example, we consider a $4-$layered medium with total length $L=7$ and its parameters are described in \autoref{table1}. The initial wave is of the form \eqref{initial_pulse} with $x_0=9,$ which is also the detector position. The measured data, meaning
the positions and heights of the major peaks (referring to single reflections), are presented in \autoref{table2}.


{\renewcommand{\arraystretch}{1.1}
\begin{table}
\parbox{.45\linewidth}{
\centering
\begin{tabular}{ccc}
\hline
Layer & Length $(\ell)$  & Wave speed $(c)$ \\\hline
1 & 2.5 & 3/7 \\
2 &  1.5 & 2/5 \\
3 & 1 & 1/2 \\
4 & 2 & 4/5 \\
\end{tabular}
\caption{Material parameters.}\label{table1}
}
\hfill
\parbox{.45\linewidth}{
\centering
\begin{tabular}{ccc}
\hline
Peak & Time $(t)$ & Height $(|h|)$ \\
\hline
1  & 4 & 0.2000 \\
2 & 47/3 & 0.0145\\
3 & 139/6 & 0.0461\\
4 & 163/6 & 0.0956\\
5 & 193/6 & 0.3923
\end{tabular}
\caption{Phaseless-data.}\label{table2}
}
\end{table}

\begin{table}
\centering
\begin{tabular}{cccccc}
\hline
Sequence & $c_1$ & $c_2$  & $c_3$ & $c_4$ &  $\sum_{j=1}^4 (t_{j+1}-t_j) c_j$ \\ \hline
1 & 3/7 & 0.400 & 0.320 & 0.200  & 10.280\\
2 & 3/7 & 0.400	 & 0.320	& 0.512 & 11.840 \\
3 & 3/7 & 0.400	& 0.500	& 0.313 & 11.562 \\
\rowcolor{Gray}
4 & 3/7 & 0.400	 & 0.500	 & 0.800 & 
14.000 \\
5 & 3/7 & 0.459 &	0.367 & 	0.229 & 11.061 \\
6 & 3/7 & 0.459	& 0.367 & 	0.588 & 12.852 \\
7 & 3/7 & 0.459 & 	0.574 &	0.359 & 12.533 \\
8 & 3/7 & 0.459 & 	0.574 &	0.918&  15.331
\end{tabular}
\caption{The output of the iterative scheme for the $4-$layered medium. All sequences solve \eqref{final_system} but only the fourth one (highlighted) satisfies in addition \eqref{eq_L}, since $2L=14.$ }\label{table3}
\end{table}

We focus at the solution of the inverse problem from phaseless data and we present the output of the proposed iterative scheme. The  system of equations takes the form:
\begin{align*}
\frac{c_1 -1}{2(c_1 +1)} &= h_1, \\
\frac{2 c_1 (c_2 - c_1)}{(c_1 +1)^2 (c_2 +c_1)} &= \pm h_2, \\
 \frac{8c_1^2 c_2 (c_3 - c_2)}{(c_2 +c_1)^2 (c_1 + 1)^2 (c_3 + c_2)} &= \pm h_3, \\
 \frac{32c_1^2 c_2^2 c_3 (c_4 - c_3)}{(c_2 +c_1)^2 (c_1 + 1)^2 (c_3 + c_2)^2 (c_4 +c_3)} &= \pm h_4
\end{align*}
In \autoref{table3} we present the $2^3$ sequences of possible solutions, one of them satisfies the condition \eqref{eq_L} and is the output of the algorithm.  Furthermore, as expected, all sequences satisfy the constraint induced by the reflection of the wave from the boundary at $x=0$, which 
reads as $\rho_5(c_1,c_2,c_3,c_4)=h_5$; thus it does not provide any extra information. 

\section{Conclusions and future work}

In this work we presented the solution of the direct and inverse scattering problem, associated with the 1-D wave equation for a multi-layered medium with constant refractive indices. To our knowledge this is the first time that the Fokas method is employed for solving this problem, producing a d'Alembert-type solution for the  direct problem on the single-layered medium. One could wonder whether this methodology could be generalised to a medium with variable refractive index (in the current work is piecewise constant), as well as to the wave equation in more spatial dimensions. For the latter question we speculate a positive answer, taking into account other works on evolution equations which apply the Fokas method to more spatial variables \cite{KF10,F02}. For the case of variable refractive index, the answer seems more complicated, but still doable taking into account the work of Deconinck et al on the heat equation \cite{DPS14,FD23}.

In this work we provided some evidence that the knowledge of the total length of the medium $L$ is enough to uniquely reconstruct the lengths and the refractive indices of the layers of medium, but it is not guaranteed that we will obtain a unique reconstruction in general for every $N-$layered medium. We note that we have no sign of failure of this strategy, but a rigorous proof is still an open question. In addition, this is an exact reconstruction method and the case of real data (data with noise) needs special treatment. The set of equations has to be replaced by a constrained minimization problem and then error and convergence analysis is needed. Both topics are out of scope of this work but an important task for future research.

\bibliographystyle{siam}
\bibliography{KalMin23}

\begin{appendices}

\section{}\label{appendix}

We present a different way to solve \eqref{cont_eq00} by splitting the time in specific sub-intervals. We restrict ourselves in the time interval $t\in (0,\, 2L/c_1 ),$ then $Q$ is given by \eqref{eq_qq0} 
and the system of equations \eqref{cont_eq00} takes the form (by integrating \eqref{cont_eq03})
\begin{subequations}\label{cont_eqb}
\begin{alignat}{3}
G (t) &= H ( t),  \quad && t>0, \label{cont_eq1b}\\
 c_1  U_0 (L+t) - c_1 G (t) &= H ( t), \quad && t>0. \label{cont_eq2b}
\end{alignat}
\end{subequations}
The solution of \eqref{cont_eqb} is given by
\begin{equation}\label{eq_g0}
G(t) = \frac{c_1}{1+c_1}  U_0 (L+t), \quad t\in \left(0, \frac{2L}{c_1} \right).
\end{equation}
 
 
As before, in the next time interval $ (2L/c_1 ,\, 4L/c_1 ),$ the system of equations \eqref{cont_eq00} takes the form
\begin{subequations}\label{cont_eqc}
\begin{alignat}{3}
G (t) &= H (t), \quad && \tfrac{2L}{c_1} \leq t< \tfrac{4L}{c_1},  \label{cont_eq0c}\\
 c_1  U_0 (L+t) - c_1 G (t) &= H ( t) + 2 H \left(t-\tfrac{2L}{c_1}\right),\quad && \tfrac{2L}{c_1} \leq t< \tfrac{4L}{c_1}. \label{cont_eq1c}
\end{alignat}
\end{subequations}
The last term in \eqref{cont_eq1c}, given \eqref{cont_eqb}, takes the form
\begin{equation}
H \left(t-\frac{2L}{c_1}\right) = 
  G \left(t-\frac{2L}{c_1}\right), \quad t > \frac{2L}{c_1}.
\end{equation}
Combining the last two equations together with \eqref{eq_g0}
results in
\begin{equation}\label{eq_g0_t2}
G (t) = 
  \frac{c_1}{1+c_1} U_0 (L+t) - 2 \frac{c_1}{(1+c_1)^2} U_0 \left( L+t-\frac{2L}{c_1} \right), \quad \frac{2L}{c_1} < t< \frac{4L}{c_1}.
\end{equation}

The presented steps can be continued to derive the boundary function $G$ for all the following time intervals. For example, in the next time interval, we get
\begin{equation}\label{eq_g0_t3b}
\begin{aligned}
G (t) &=   \frac{c_1}{1+c_1} U_0 (L+t) - 2 \frac{c_1}{(1+c_1)^2}U_0 \left(L+t-\frac{2L}{c_1}\right),\\
&\phantom{=}+ 2 \frac{c_1 (1-c_1)}{(1+c_1)^3}U_0 \left(L+t-\frac{4L}{c_1}\right), \quad \frac{4L}{c_1} < t< \frac{6L}{c_1}.
\end{aligned}
\end{equation}
 
 This process can be generalized to obtain the function $G$ for all time intervals, leading to the form:
 \begin{equation}\label{g0_final}
\begin{aligned}
G(t) &= 
 \frac{c_1}{1+c_1} U_0 (L+t) - 2 c_1 \sum_{n=0}^{N} \frac{(c_1-1)^n}{(c_1+1)^{n+2}}U_0 \left(L+t-\frac{2(n+1)L}{c_1}\right), \\
 &\phantom{=}\mbox{for}
 \quad \frac{2(N+1)L}{c_1}<t  < \frac{2(N+2)L}{c_1}, \quad N=-1,0,1,\ldots,
\end{aligned}
\end{equation}
 with $N=-1,$ we mean that no summation is performed.
 
 With similar but more lengthy calculations we obtain the function $G_0$ also for a double-layered medium.



\section{}\label{appendix-b}

We consider \eqref{layer2_final1} once for $t = D+ \tfrac{2l_1}{c_1}+ \tfrac{2l_2}{c_2},$ and then for $t = D+ \tfrac{2l_2}{c_2},$ and we subtract the two resulted formulas to derive
\begin{equation}\label{eq_g0_2l1_2l2}
G_0^{(2,2)} - G_0^{(0,2)} = -\frac{2 }{c_1 +1} G_0^{(0,2)} 
+ \frac{2 }{c_1 +1} G_1^{(1,2)}.
\end{equation}

In order to avoid technicalities in the presentation of the computation of the terms involved in the above equation, we make the  assumption that the medium is such that $1<\frac{c_1 l_2}{c_2 l_1} <2.$ 
From a physical point of view, this constraint means that the major peaks will arrive before the minor peaks.
We emphasise that the other cases can be treated similarly, yielding identical results.

The 
Equation \eqref{layer2_final1} gives
\begin{equation}\label{eq_g0_l2}
G_0^{(0,2)} = -\frac{2}{c_1 +1} G_0^{(-2,2)} +\frac{2}{c_1 +1} G_1^{(-1,2)},
\end{equation}
since the remaining remaining terms are zero from \eqref{conditions}.  In the right-hand side, the first term is zero. The second term is given by \eqref{eq_2l_g1} where we observe that the leading terms are again zero, thus also the second term in \eqref{eq_g0_l2} is zero, resulting in
\begin{equation}\label{eq_g0_l2b}
G_0^{(0,2)} = 0.
\end{equation}
The remaining term in \eqref{eq_g0_2l1_2l2}, using \eqref{eq_2l_g1} and taking into account the above assumptions, admits the form
\begin{equation}\label{eq_g1_l2}
\begin{aligned}
G_1^{(1,2)} &= - \frac{2 c_1}{c_2 +c_1} G_1^{(1,0)} \\
&\overset{\eqref{eq_2l_g1}}{=} - \frac{2 c_1}{c_2 +c_1} \frac{2 c_2}{c_2 +c_1} G_0^{(0,0)} \\
&\overset{\eqref{eq_2l_g0}}{=} - \frac{4c_1^2 c_2}{(c_2 +c_1)^2 (c_1 + 1)}.  
\end{aligned}
\end{equation}
Finally, from \eqref{eq_g0_2l1_2l2} considering 
\eqref{eq_g0_l2b} and \eqref{eq_g1_l2}, 
we obtain  \eqref{eq_g0_l1_l2}.

\end{appendices}

\end{document}